\documentclass[11pt]{article}
\usepackage{latexsym}

\newcommand{\be}{\begin{equation}}
\newcommand{\ee}{\end{equation}}
\newcommand{\R}{\rm I\kern-.19emR}

\newcommand{\bd}{\begin{displaymath}}
\newcommand{\ed}{\end{displaymath}}

\newtheorem{1}{LEMMA}[section]

\newtheorem{2}[1]{THEOREM}

\newtheorem{6}[1]{DEFINITION}

\newtheorem{10}[1]{COROLLARY}

\newcommand{\bp}{\underline{\bf Proof}:\ }
\newcommand{\ep}{{\hfill $\Box$}\\ }
\def\N{{\rm I\kern-.25em N}}

  \def\C{{\rm\kern.24em
\vrule width.02em height1.4ex depth-.05ex \kern-.26em C}}

\begin{document}

\title{Hermitian Positive Semidefinite Matrices Whose Entries
Are $0$ Or
$1$ in Modulus\footnote{Research supported in part by NSF Grant
DMS-942436}}
 \author{  Daniel Hershkowitz \\
           Department of Mathematics \\
           Technion -- Israel Institute of Technology \\
            Haifa  32000, Israel
   \and Michael Neumann \\
        Department of Mathematics \\ University of Connecticut \\
        Storrs, Connecticut  06269--3009, USA
    \and Hans Schneider \\
     Department of Mathematics \\
     University of Wisconsin  \\
     Madison, Wisconsin 53706, USA }
\date{10 July 1998}
\maketitle

\begin{abstract}
We show that a matrix is a  Hermitian positive
semidefinite matrix
whose nonzero entries have modulus $1$ if and only if it similar
to a
direct sum of all $1's$ matrices and a $0$ matrix via a
unitary monomial similarity.  In
particular,
the only such  nonsingular matrix  is the
identity matrix
and the only such irreducible matrix is similar to an all
$1$'s matrix by means of a unitary diagonal similarity.
Our results extend earlier results of Jain and Snyder for the
case in which
the nonzero entries (actually)  equal $1$. Our methods of proof,
which rely on the so called  principal submatrix rank
property,
differ
from the approach used by Jain and Snyder.
\end{abstract}

\thispagestyle{empty}

\newpage

\section{Introduction}
\label{sec.1}

In this note we characterize the set of all Hermitian
positive semidefinite matrices $A$
whose entries have  modulus $1$ or
$0$. We comment that the
real matrices in this set whose diagonal
entries are all nonzero, and hence necessarily $1$, belong
to the set of {\bf correlation matrices} (see Horn and Johnson
\cite[p.400]{HJ}).\\

Our characterization here is motivated by a surprising result due
to
Jain and Snyder \cite{JS} which can be stated as follows:\\

\begin{2}
\label{thm.js}
{\rm (Jain and Snyder \cite[Theorems 2 and 3]{JS})}
If $A$ is any positive semidefinite $(0,1)$--matrix, then $A$ is
permutationally
similar to a direct sum of matrices each of which is either an
all $1$'s matrix or a
zero matrix.
In particular,
if $A$ is {\rm (}also{\rm )} irreducible, then $A$ is the
$n\times n$ all
$1$'s matrix.
\end{2}

Jain's and Snyder's proof of Theorem \ref{thm.js} rests on
their observation
that any
positive semidefinite $(0,1)$--matrix has a Cholesky
factorization with
a $(0,1)$--Cholesky factor. The
proof of our generalization of Theorem \ref{thm.js} relies on the
following property which is possessed by Hermitian positive
semidefinite
matrices and from which results on Cholesky factorizations
follow:\\

\begin{6}
\label{def.psrp}
 {\rm A matrix $A\in \C^{n,n}$ is said to have the
{\it principal submatrix rank property} {\rm (}PSRP{\rm )}
if the
following conditions holds:\\

(i) The column space determined by every set of rows of $A$ is
equal
to the column space of the principal submatrix lying in these
rows.\\

(ii) The row space determined by every set of columns of $A$ is
equal
to the row space of the principal submatrix lying in these
columns.}
\end{6}

It is known (cf. Hershkowitz and Schneider \cite{HS,HS1}
and Johnson \cite{J})
that a positive semidefinite
Hermitian matrix has PSRP. We give the following simple proof for
the sake of completeness:\\

Since a
permutation similarity applied to a positive semidefinite
Hermitian matrix
again yields a positive semidefinite Hermitian matrix,
it is enough to show that the row space
determined by the first
$k$ columns of a positive semidefinite Hermitian matrix $A$ is
equal to the
row space determined by the leading principal minor $B$ of $A$ of
size $k$. This is equivalent to showing the following: If $v^T =
[w,0]^T$, where $w$ is of length $k$, and $Bw = 0$ then $Av = 0$.
But, for such $v$, we have $v^*Av = w^*Bw = 0$, and the result
follows since $A$ is positive semidefinite.\\

\section{Main Result}
\label{sec.2}

To facilitate our extension
of Theorem \ref{thm.js}, we
introduce the following notion:

\begin{6}
\label{def.1}
{\rm A matrix $P\in \C^{n,n}$ is called a {\em unitary monomial
matrix} if $P = QD$, where $Q$ is a permutation matrix and $D$ is
a diagonal matrix all of whose diagonal entries are of modulus
$1$.}
\end{6}

We are now ready to state the main result of this note:\\

\begin{2}
\label{diningroom}
 A matrix $A \in \C^{n,n}$ is Hermitian  positive
semidefinite and all its entries  have modulus $1$ or $0$
if and only if  $A$ is  similar, by means of a unitary monomial
matrix,  to a direct sum of matrices each of which is either an
all $1$'s matrix or a
zero matrix.
\end{2}

\bp
The proof of the ``if'' part is obvious, so we proceed to prove
the ``only if'' part.
This is done by induction on $n$, the size
of $A$. The result is trivial if $n=1$. So let $n > 1$ and assume
theorem  holds for matrices of all sizes less than $n$.\\

If $A$ is reducible, then $A$ is permutation similar to a direct
sum of
positive semidefinite Hermitian matrices of size less than $n$
whose
nonzero entries are all of modulus $1$. By our inductive
assumption each
direct summand is unitarily monomially similar to the direct sum
of all $1$'s
matrices and a $0$ matrix, and hence the same is true for $A$.\\

So assume that $A$ is irreducible. We shall first show that there
exists
$(0,1)$--matrix $E = D^{-1}P^{-1}APD$, where $D$ is diagonal and
$P=(p_{i,j})$ is a
unitary monomial matrix with $p_{n,n} = 1$. Further all elements
of the
last row and column of $E$ are $1$.\\

No diagonal element of $A$ is $0$, for then the corresponding row
and
column would also be $0$.  Since the leading submatrix $B$ of
size $n-1$
of $A$ is positive semidefinite, it follows there is a unitary
monomial
similarity of $B$ such that the resulting matrix is a direct sum
of all
$1$'s matrices and a $0$ matrix. We extend this similarity to a
unitary
monomial similarity of $A$ which leaves the last row and column
of $A$ in
place. We thus obtain a matrix $C = P^{-1}AP$, where $p_{n,n} =
1$, such
that all nonzero elements in the last row and column of $C$ are
of modulus
1 and the leading submatrix of size $n-1$ of $C$ is a direct sum
of all
$1$'s matrices and a $0$ matrix. We partition the last row and
column of
$C$ in conformity with this direct sum. Since $C$ is positive
semidefinite, it follows by PSRP that each subvector of the last
row and
column determined by this partition is a multiple of the all
$1$'s vector
of the appropriate size by a number of modulus $1$ or by $0$.
But if one of the last column is a $0$ multiple of the all $1$'s
vector,
then is $C$ reducible. Hence each subvector of the last column is
a
multiple of an all $1$'s vector by a number of modulus $1$. Let
$D$ be the
unitary diagonal matrix whose diagonal entries coincide with the
last
column of $C$. Since $D$ has equal entries corresponding to the
blocks of
$B$, it follows that $E = D^{-1}CD$ is a $(0,1)$--matrix whose
last row
and column consists of $1$'s. \\

Our proof shows that the last row and column of $A$ have no $0$
entries. By applying permutation similarities  to $A$ and
repeating the above construction, we deduce that the same is true
of every row and column. Hence the matrix $E$ we have obtained is
the all $1$'s matrix.
\ep

Theorem \ref{diningroom} has several corollaries:\\

\begin{10}
\label{cor.1}
A matrix $A \in \C^{n,n}$ is a positive
definite Hermitian matrix whose nonzero entries are all of
modulus $1$ if and only if
$A$ is the identity matrix.
\end{10}

\begin{10}
\label{cor.2}
A matrix  $A \in \C^{n,n}$ is an irreducible positive
semidefinite Hermitian matrix whose
nonzero entries are all of
modulus $1$ if and only if $A$ can be transformed
to the all $1$'s matrix by a unitary diagonal similarity.
\end{10}

\begin{10}
\label{cor.3}
Let $A \in \C^{n,n}$ be a positive
semidefinite Hermitian matrix whose nonzero entries are all of
modulus $1$. Then there is an LU factorization of $A$ with $L$
nonsingular where $L$
and $U$ are similar to $(0,1)$--matrices via the same unitary
diagonal matrix.
\end{10}

\bp If $C$ is the $k \times k$ block of all $1$'s
then
it admits the factorization $C = L_1U_1$, where the first column
of
$L_1$ consists of $1$'s, the diagonal entries of $L_1$ are $1$,
the
first row of $U_1$ consists of $1$'s and all other elements of
$L_1$
and $U_1$ are $0$. Now, if $C$ is a direct sum of such blocks and
a $0$
matrix, it easily follows that $C$ admits an LU factorization
where $L_2$ is a lower triangular nonsingular $(0,1)$--matrix
and $U_2$ is an upper triangular $(0,1)$--matrix. If $B$ is
permutation similar to $C$ then we can find a permutation matrix
$P$ which does not change the order of the order of rows  and
columns in any given block and for which $B = P^TCP$. Then
$L_3 = P^TL_2P$ is a lower triangular nonsingular $(0,1)$--matrix
and $U_3 = P^TU_2P$ is an upper triangular $(0,1)$--matrix. If $A
= D^*BD$, where $D$ is a unitary diagonal matrix, it follows that
$L = D^*L_3D$ and $U = D^*U_3D$ satisfy the conditions of the
corollary. The conclusion of the corollary now follow by applying
Theorem \ref{diningroom}.
\ep

In a very similar way to the proof of the above corollary, we
can prove the following corollary:\\

\begin{10}
\label{cor.4}
Let $A \in \C^{n,n}$ be a positive
semidefinite Hermitian matrix whose nonzero entries are all of
modulus $1$. Then there is an Cholesky $LL^*$ factorization of
$A$ where $L$ is similar to $(0,1)$--matrices via a unitary
diagonal matrix.
\end{10}

\vspace{.3in}
\hspace{-.25in} \underline{\bf ACKNOWLEDGEMENT} \  The authors
wish to thank
Professor
Emeric Deutsch for bringing the question of characterizing the
$n\times n$ $(0,1)$--matrices which are positive semidefinite to
their attention.

\end{document}